
\input epsf.tex
\input amssym.def
\input amssym
\magnification=1100
\baselineskip = 0.22truein
\lineskiplimit = 0.01truein
\lineskip = 0.01truein
\vsize = 8.7truein
\voffset = 0.2truein
\parskip = 0.10truein
\parindent = 0.3truein
\settabs 12 \columns
\hsize = 6.0truein
\hoffset = 0.1truein

\setbox\strutbox=\hbox{%
\vrule height .708\baselineskip
depth .292\baselineskip
width 0pt}
\font\caps=cmcsc10
\font\bigtenrm=cmr10 at 14pt

\def\sqr#1#2{{\vcenter{\vbox{\hrule height.#2pt
\hbox{\vrule width.#2pt height#1pt \kern#1pt
\vrule width.#2pt}
\hrule height.#2pt}}}}
\def\square{\mathchoice\sqr46\sqr46\sqr{3.1}6\sqr{2.3}4}

\centerline{\bigtenrm THE CROSSING NUMBER OF SATELLITE KNOTS}
\tenrm
\vskip 14pt
\centerline{MARC LACKENBY}
\vskip 18pt
\centerline{\caps 1. Introduction}
\vskip 6pt

One of the most basic invariants of a knot $K$ is its crossing number $c(K)$,
which is the minimal number of crossings in any of its diagrams. However, it remains
quite poorly understood. For example, it is a notorious unsolved conjecture that
if $K_1 \sharp K_2$ is the connected sum of two knots $K_1$ and $K_2$, then
$c(K_1 \sharp K_2) = c(K_1) + c(K_2)$. Connected sums are particular cases
of satellite knots, which are defined as follows. Let $L$ be a non-trivial knot in the 3-sphere. 
Then a knot $K$ is a {\sl satellite knot} with {\sl companion knot} $L$
if $K$ lies in a regular neighbourhood $N(L)$ of $L$, it does not
lie in a 3-ball in $N(L)$ and is not a core curve of $N(L)$. (See Figure 1.) It is conjectured
that the $c(K) \geq c(L)$ (Problem 1.67 in Kirby's problem list [3]). In this paper, 
we establish that an inequality of this form holds, up to a universally bounded factor.

\noindent {\bf Theorem 1.1.} {\sl Let $K$ be a satellite knot with companion knot $L$. Then
$c(K) \geq c(L)/10^{13}$.}

\vskip 18pt
\centerline{
\epsfxsize=3.5in
\epsfbox{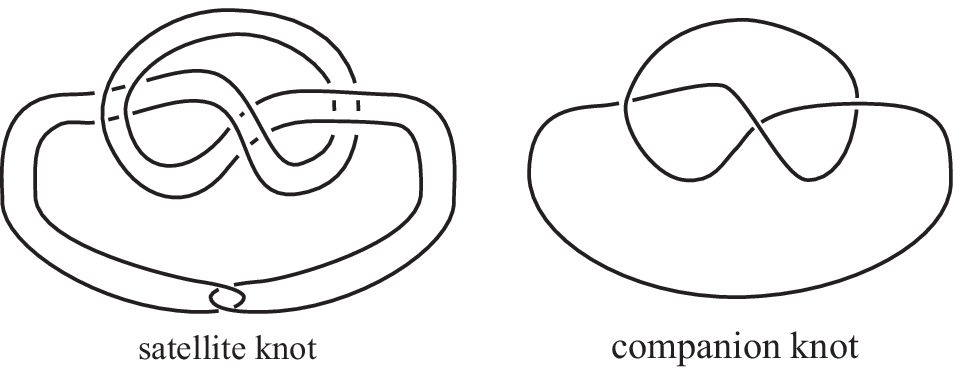}
}
\centerline{Figure 1.}

This result should be compared with the main theorem of [4], which is as follows.

\noindent {\bf Theorem 1.2.} {\sl Let $K_1 \sharp \dots \sharp K_n$ be the connected
sum of knots $K_1, \dots, K_n$. Then
$${c(K_1) + \dots + c(K_n) \over 152} \leq c(K_1 \sharp \dots \sharp K_n) \leq c(K_1) + \dots + c(K_n).$$}

Both the statement and the proof of Theorem 1.2 will be crucial for us in this
paper. Aside from this result, there has not been much work before on the crossing
number of satellite knots. In [2], Freedman and He defined the `asymptotic
crossing number' of a knot $L$ in terms of the crossing number of certain satellites
of $L$. They related this invariant to a physically defined quantity called
the `energy' of $L$, and also showed that the asymptotic crossing number is bounded
below by a linear function of the genus of $L$. They conjectured that the asymptotic
crossing number is equal to the crossing number, a result which would follow from
various stronger versions of Kirby's Problem 1.67. It is conceivable that the
methods behind the proof of Theorem 1.1 may be applied to obtain new information
about the asymptotic crossing number of a knot.

Another key input in the proof of Theorem 1.1 is the machinery developed by the author in [5]. 
The main goal of that paper was to show that, given any triangulation of the solid torus, there is a core
curve (or, more precisely, a `pre-core curve') that lies in the 2-skeleton and
that intersects the interior of each face in at most $10$ straight arcs. However, in this paper,
it is more convenient to use handle structures. In order to be able
to speak of `straight' arcs in a handle decomposition, we give it an
`affine' structure, which is defined as follows.

Whenever we refer to a handle structure on a 3-manifold, we insist that each
handle is attached to handles of strictly lower index.
An {\sl affine handle structure} on a 3-manifold $M$ is a handle structure where each 
0-handle and 1-handle is identified with a compact (but possibly non-convex) polyhedron in ${\Bbb R}^3$, 
so that 
\item{(i)} each face of each polyhedron is convex;
\item{(ii)} whenever a 0-handle and 1-handle intersect, each component of intersection is
identified with a convex polygon in ${\Bbb R}^2$, in such a way that the inclusion
of this intersection into each handle is an affine map with image equal to a face of the
relevant polyhedron;
\item{(iii)} for each 0-handle $H_0$, each component of intersection with a 2-handle,
3-handle or $\partial M$ is a union of faces of the polyhedron associated with $H_0$;
\item{(iv)} the polyhedral structure on each 1-handle is the product of a convex 2-dimensional
polygon $P$ and an interval $I$, where $P \times \partial I$ is the intersection with the
0-handles, and the intersection between the 1-handle and any 2-handle is $\beta \times I$,
where $\beta$ is a union of disjoint sides of $P$.

\noindent Since each 0-handle and 1-handle is identified with a polyhedron,
it makes sense to speak of a straight arc in that handle.

The following result was Theorem 4.2 of [5].

\noindent {\bf Theorem 1.3.} {\sl Let ${\cal H}$ be an affine handle structure of the
solid torus $M$. Suppose that each $0$-handle of ${\cal H}$ has at most $4$ components
of intersection with the $1$-handles, and that each $1$-handle has at most $3$
components of intersection with the $2$-handles. Then $M$ has a core curve that intersects 
only the $0$-handles and $1$-handles, that respects the product structure on the
$1$-handles, that intersects each $1$-handle in at most $24$ straight arcs,
and that intersects each $0$-handle in at most $48$ arcs. Moreover, the
arcs in each $0$-handle are simultaneously parallel to a collection of arcs $\alpha$ in the boundary of the
corresponding polyhedron, and each component of $\alpha$ intersects each face of the polyhedron
in at most $6$ straight arcs.}

The proof of Theorem 1.1 follows a similar route to that of Theorem 1.2, but it requires some
new ideas. We start with a diagram $D$ for the satellite knot $K$ with minimal crossing number.
It would suffice to construct a diagram for the companion $L$ with crossing number at most
$10^{13} \ c(D)$. However, instead, we construct a diagram $D'$
for a knot $L'$ which has $L$ as a connected summand, and with
crossing number at most $3 \times 10^{10} c(D)$. Then, applying Theorem 1.2,
we deduce that
$$c(L) \leq 152 \ c(L') \leq 152 \ c(D') \leq 152 \times 3 \times 10^{10} \ c(D) < 10^{13} \  c(K),$$
as required.

This diagram $D'$ is constructed as follows.
We use the diagram $D$ to build a handle structure ${\cal H}_X$ for the exterior $X$ of $K$.
Let $T$ be the torus $\partial N(L)$ arising from the satellite construction.
Then $T$ is essential in $X$, and so may be placed in normal form with respect
to ${\cal H}_X$. Now cut $X$ along $T$, to give two 3-manifolds, one
of which is a copy of the exterior of $L$, the other of which is denoted
$Y$ and is $N(L) - {\rm int}(N(K))$. The aim is to find a handle structure
${\cal H}'_{Y'}$ for $Y$ which sits nicely in ${\cal H}_X$. In particular,
the 0-handles (respectively, 1-handles) of ${\cal H}'_{Y'}$ lie in the
0-handles (respectively, 1-handles) of ${\cal H}_X$, and each 0-handle
and 1-handle of ${\cal H}_X$ contains at most $6$ such handles of ${\cal H}'_{Y'}$. Now,
$Y$ inherits a handle structure ${\cal H}_Y$ from ${\cal H}_X$, but it
may not have the required properties. This is because the torus $T$ may intersect
a handle $H$ of ${\cal H}_X$ in many normal discs, and these may divide $H$
into many handles of ${\cal H}_Y$. However, the normal discs come in
only finitely many types (at most $5$, in fact, in any handle), and normal discs of
the same type are parallel. Between adjacent parallel normal discs,
there is an $I$-bundle, and these patch together to form the {\sl parallelity
bundle} for ${\cal H}_Y$. This structure was considered in detail in [4],
where it was shown that this may be enlarged to another $I$-bundle ${\cal B}$
called a {\sl generalised parallelity bundle}. This has many nice properties.
In particular, one may (under certain circumstances) ensure that it consists of
$I$-bundles over discs, and other components which have incompressible
vertical boundary, which are annuli properly
embedded in $Y$, with boundary in $T$. Now, it would be convenient if there
were no embedded essential annuli in $Y$ with boundary in $T$, but unfortunately
there may be. This happens, for example, if $K$ is also a satellite of
a knot $L'$ which has $L$ as a connected summand. (An example is shown in Figure 10.)
We hypothesise this situation
away, by focusing instead on $L'$. This is why we aim to find a diagram
for $L'$ instead of $L$. Using this line of argument, and others,
we arrange that $Y$ contains no properly embedded essential annuli
with boundary in $T$. Hence, ${\cal B}$ consists of $I$-bundles over discs.
We replace each of these by a 2-handle, thereby constructing the required handle structure ${\cal H}'_{Y'}$.
We now attach the solid torus $N(K)$, forming a handle structure ${\cal H}_{V'}$
for the solid torus $N(L')$. Then, using Theorem 1.3, we find a core curve for $N(L')$
which lies nicely with respect to ${\cal H}_{V'}$ and hence ${\cal H}_X$.
This is a copy of $L'$, and projecting, we obtain the diagram $D'$ for $L'$ with
the required bound on crossing number.

The factor $10^{13}$ is very large, and one may wonder whether there are
ways of reducing it. We have not attempted to optimise this constant. In general,
where there was a choice between two arguments, one shorter and simpler than the
other, but with worse constants, we have opted for the short and simple route.
However, we would be surprised if the constant could be reduced by more
than a factor of $10^8$ without a major modification to the argument.
Nevertheless, there are versions of the theorem with significantly
improved constants, but with weakened conclusions. These are as follows.

\noindent {\bf Theorem 1.4.} {\sl Let $K$ be a satellite knot with companion knot
$L$. Suppose that there is no essential torus properly embedded in $N(L) - {\rm int}(N(K))$.
Then,
$$c(K) \geq {c(L) \over 3 \times 10^{10}}.$$}

\noindent {\bf Theorem 1.5.} {\sl Let $K$ be a satellite knot with companion knot $L$.
Suppose that $L$ is prime. Then, for some knot $\tilde L$, which is either $L$ or a cable of $L$,
$$c(K) \geq {c(\tilde L) \over 152}.$$}

The plan of this paper is as follows. In Section 2, we explain how a handle structure
on the exterior of a link $K$ can be constructed from any connected diagram for $K$.
This is similar to a construction in [4]. However, we then place an affine structure
upon this handle decomposition, which is new. In Section 3, we recall some of the
normal surface that we will need, including the notion of a parallelity bundle
and generalised parallelity bundle from [4]. In Section 4, we give the proof of
the main theorem, but with one step excluded. It turns out that before the main
theorem can be proved in full generality, the special case of certain 2-cables
must be analysed separately. We do this in Section 5. Finally, in Section 6, we
explain how the proof of the main theorem can be adapted to give Theorem 1.5.

\vskip 18pt
\centerline{\caps 2. An affine handle structure from a diagram}
\vskip 6pt

In [4], we introduced a method for creating a handle structure for
the exterior of a link $K$, starting from a connected diagram $D$ for $K$.
In this section, we recall the main details of this construction.

The diagram is a 4-valent graph embedded in the 2-sphere, and we
realise this 2-sphere as the equator in $S^3$.
Let $S^2 \times [-1,1]$ be a regular neighbourhood of this 2-sphere,
where $S^2 \times \{ 0 \}$ is the equator itself.

The diagram specifies an embedding of $K$ into the 3-sphere, so that
away from small neighbourhoods of the crossings, it lies in the
diagram 2-sphere, and at each crossing, two arcs of $K$ come out of the
diagram 2-sphere. One goes vertically upwards to height $1$, then
runs horizontally, and then returns to the diagram 2-sphere. The other
arc makes a similar itinerary below the diagram. Thus, $K$ lies in
$S^2 \times [-1,1]$ and its image under the product projection map
to $S^2$ equals the 4-valent graph specified by $D$.

The 0-handles and 1-handles of the handle structure are thickenings of
a graph that lies in the equatorial 2-sphere, as follows.
There are four 0-handles arranged around each crossing, as in Figure 2.
These are joined by four 1-handles which form a square that surrounds
the crossing. This square is small enough so that $K$ misses
these 1-handles, because it lies above and below the diagram
at these points. In addition, there are two 1-handles which follow
each edge of the diagram, and lie either side of that edge.
(See Figure 2.)

\vskip 18pt
\centerline{
\epsfxsize=3.5in
\epsfbox{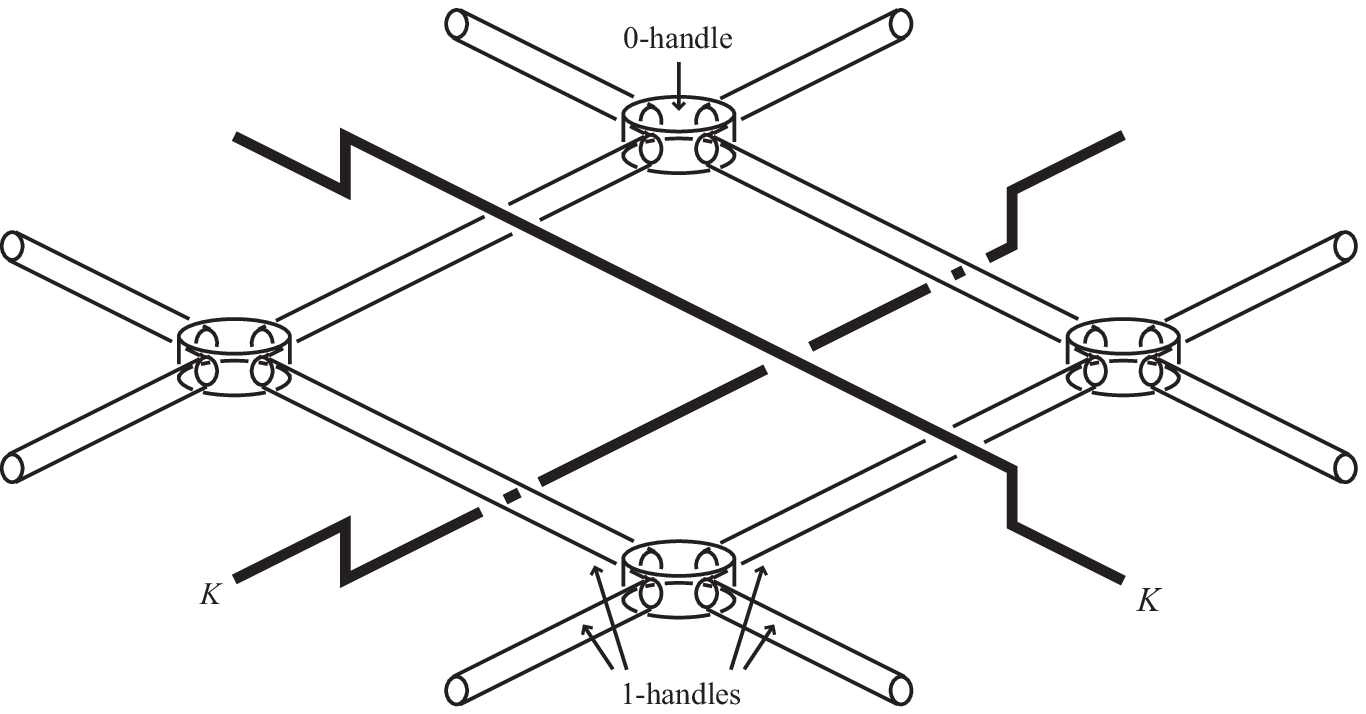}
}
\centerline{Figure 2.}

We will not describe in detail here how the 2-handles are
attached to the 0-handles and 1-handles. We merely note that they
intersect each 0-handle as shown in the left of Figure 5. 

Finally, there are two 3-handles, one of which lies entirely
above the plane of the diagram, and one of which lies below.

In [4], two modifications were made to this handle structure, and it is
convenient to make one of the same modifications here. Thus, two new 0-handles are introduced,
which lie either side of an edge of the diagram. The insertion of these
handles has the effect of dividing two 1-handles each into two.
The two new 0-handles are known as {\sl exceptional}. Two new 1-handles are
also introduced, each of which joins the two exceptional 0-handles.
These lie above and below the plane of the diagram. They lie in
2-handles of the old handle structure, and so each 1-handle
subdivides the 2-handle into two. (See Figure 3.)

\vskip 18pt
\centerline{
\epsfxsize=4.9in
\epsfbox{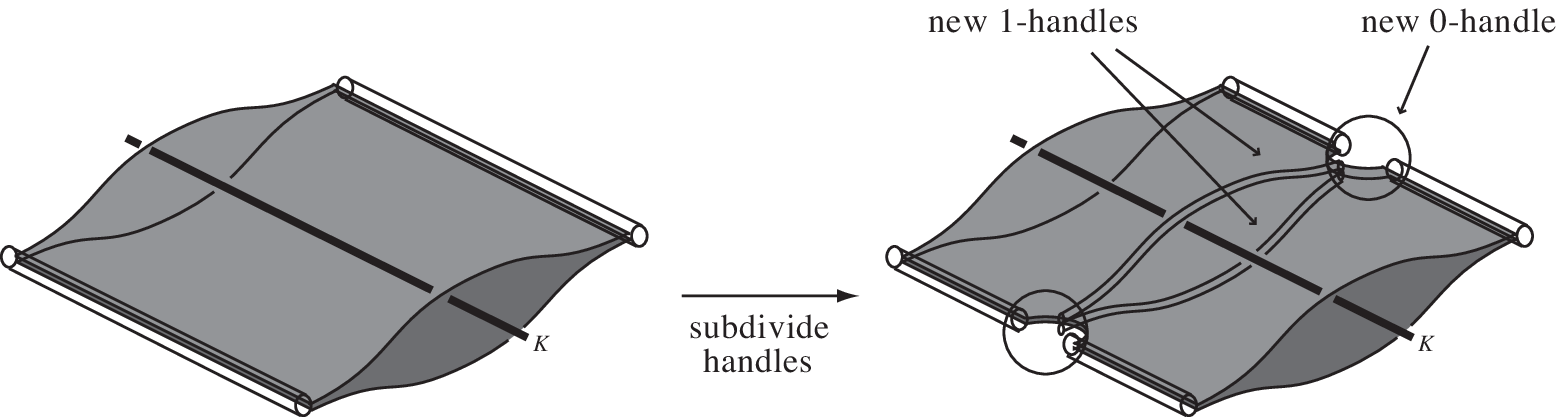}
}
\centerline{Figure 3.}

We denote this handle structure on the exterior of $K$ by ${\cal H}_X$.

Note that there is a slight discrepancy between the handle structure here
and the one considered in [4]. There, a further modification was made
which removed the two 3-handles, together with two 2-handles. We will
not take this step here.

As in [4], we want to work in ${\Bbb R}^3$ rather than $S^3$. We therefore pick a
point in the diagram 2-sphere, distant from the crossings, and declare that it is the point
at infinity. We thereby obtain a knot diagram in ${\Bbb R}^2$ which sits inside
${\Bbb R}^3$, and we may assume that the diagrammatic projection map is just the
vertical projection from ${\Bbb R}^3$ to ${\Bbb R}^2$.

We now wish to place an affine structure upon ${\cal H}_X$.
The first step is to realise each 1-handle as a polyhedron.
There are two types of 1-handle: those that form part of a square
surrounding a crossing, and those that are parallel to an edge
of the diagram. We give each of these a slightly different polyhedral
structure. Each is the metric product of a convex 2-dimensional polygon
and an interval, but the 2-dimensional polygon is a little different
in each of the two cases. The precise polygons are shown in Figure 4.

\vskip 18pt
\centerline{
\epsfxsize=3.3in
\epsfbox{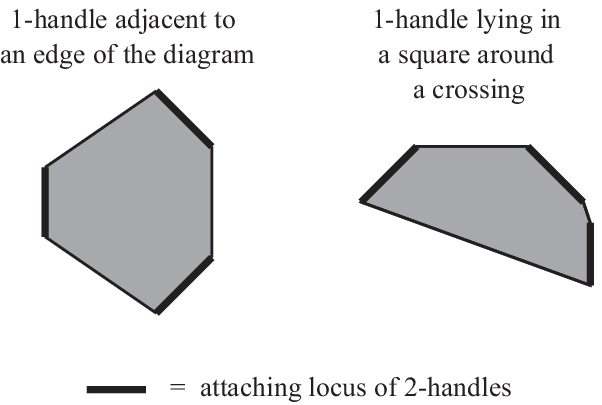}
}
\centerline{Figure 4.}
\vfill\eject

We now realise each 0-handle as a polyhedron in ${\Bbb R}^3$.
We focus on the unexceptional 0-handles. Currently, each is of the form
shown in the left of Figure 5. We replace this with the (non-convex)
polyhedron shown in the right of Figure 5. Each component of $\partial {\cal H}^0_X \cap \partial X$
is realised as a union of 4 triangular faces of the polyhedron (which are not shown in
Figure 5).

\vskip 18pt
\centerline{
\epsfxsize=3.8in
\epsfbox{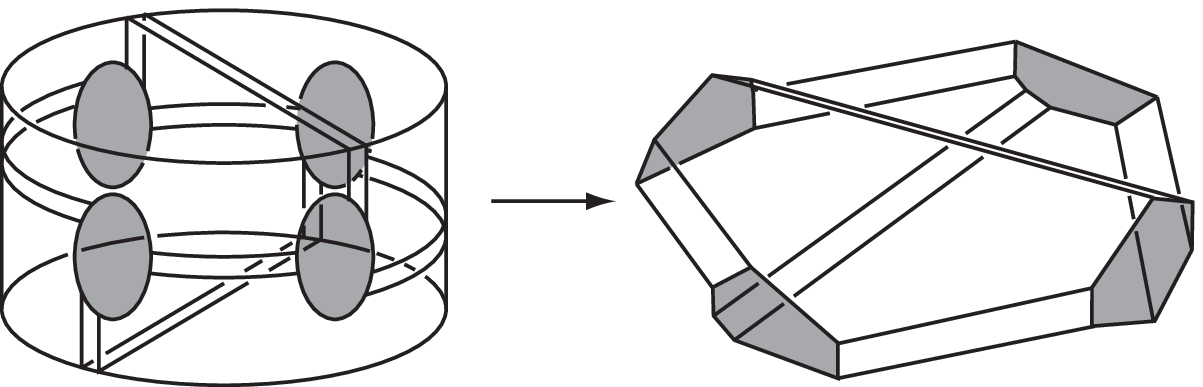}
}
\centerline{Figure 5.}

Each polyhedral 0-handle is embedded in ${\Bbb R}^3$ isometrically.
But this might not be possible for some of the 1-handles. Some 1-handles
follow an edge of the diagram and this edge might not be straight. Nevertheless,
there is a way to embed the polyhedron horizontally into ${\Bbb R}^3$,
so that it has the following property. If the 1-handle is $D^1 \times D^2$,
then for any two distinct points $x_1$ and $x_2$ in $D^2$, the arcs
$D^1 \times \{ x_1 \}$ and $D^1 \times \{ x_2 \}$ vertically
project to arcs in the diagram that either are equal or do not cross.

\vskip 18pt
\centerline{\caps 3. Generalised parallelity bundles}
\vskip 6pt

As in the proof of Theorem 1.2 that is given in [4], a key technical tool in this paper
is the notion of a generalised parallelity bundle. In this section, we will recall
this notion. We will also extend
some of the results in [4] so that they can be applied in the context of this paper.

When ${\cal H}$ is a handle structure of a 3-manifold, ${\cal H}^i$ will denote the
union of the $i$-handles.

\noindent {\bf Convention 3.1.} We will insist throughout this paper that 
any handle structure ${\cal H}$ on a 3-manifold satisfies the following conditions:
\item{(i)} each $i$-handle $D^i \times D^{3-i}$ intersects $\bigcup_{j \leq i-1} {\cal H}^j$
in $\partial D^i \times D^{3-i}$;
\item{(ii)} any two $i$-handles are disjoint;
\item{(iii)} the intersection of any 1-handle
$D^1 \times D^2$ with any 2-handle $D^2 \times D^1$ is of
the form $D^1 \times \alpha$ in $D^1 \times D^2$, where
$\alpha$ is a collection of arcs in $\partial D^2$,
and of the form $\beta \times D^1$ in $D^2 \times D^1$,
where $\beta$ is a collection of arcs in $\partial D^2$;
\item{(iv)} each 2-handle of ${\cal H}$ runs over at least one 1-handle.

\noindent The handle structure constructed in Section 2 satisfies these
requirements.

Now let $X$ be a compact orientable 3-manifold with a handle structure ${\cal H}_X$.

Let ${\cal F}$ be the surface ${\cal H}^0_X \cap ({\cal H}^1_X \cup {\cal H}^2_X)$,
let ${\cal F}^0$ be ${\cal H}^0_X \cap {\cal H}^1_X$, and 
let ${\cal F}^1$ be ${\cal H}^0_X \cap {\cal H}^2_X$.
By the above conditions, ${\cal F}$ is a thickened graph,
where the thickened vertices are ${\cal F}^0$ and the thickened
edges are ${\cal F}^1$.

\noindent {\bf Definition 3.2.} 
We say that a surface $T$ properly embedded in $X$ is {\sl standard}
if
\item{(i)} it intersects each 0-handle in a collection of properly
embedded disjoint discs;
\item{(ii)} it intersects each 1-handle $D^1 \times D^2$ in 
$D^1 \times \beta$, where $\beta$ is a collection of properly
embedded disjoint arcs in $D^2$;
\item{(iii)} it intersects each 2-handle $D^2 \times D^1$ in 
$D^2 \times P$, where $P$ is a collection of points in the
interior of $D^1$;
\item{(iv)} it is disjoint from the 3-handles.

\noindent See Figure 6.

\vskip 18pt
\centerline{
\epsfxsize=4in
\epsfbox{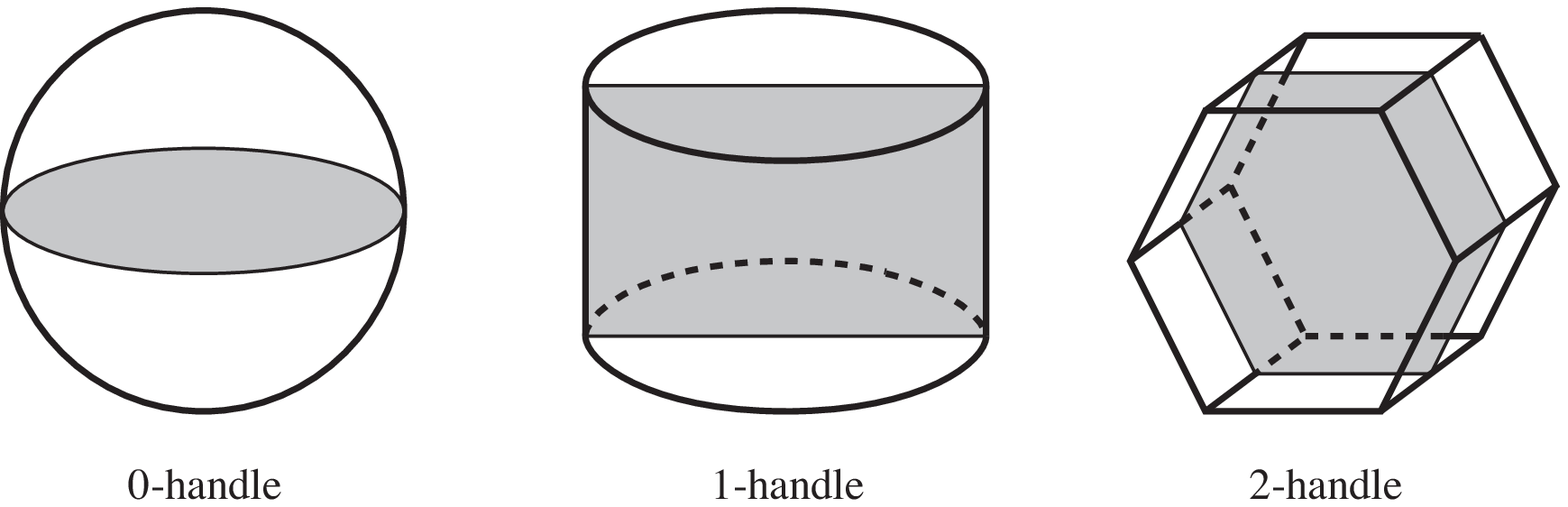}
}
\vskip 6pt
\centerline{Figure 6.}

\noindent {\bf Definition 3.3.} A disc component $D$ of $T \cap {\cal H}^0_X$ is
said to be {\sl normal} if
\item{(i)} $\partial D$ intersects any thickened edge of ${\cal F}$ in at most one arc;
\item{(ii)} $\partial D$ intersects any component of $\partial {\cal F}^0 - {\cal F}^1$ at most once;
\item{(iii)} $\partial D$ intersects each component of $\partial {\cal H}^0 - {\cal F}$ in at most
one arc and no simple closed curves.

\noindent A standard surface that intersects each 0-handle in a disjoint union of normal
discs is said to be {\sl normal}. (See Figure 7.)

This is a slightly weaker definition of normality than is used by some authors,
for example Definition 3.4.1 in [6]. However, if we had used the definition in [6],
Proposition 3.4 (below) would no longer have held true.

\vskip 12pt
\centerline{
\epsfxsize=2.6in
\epsfbox{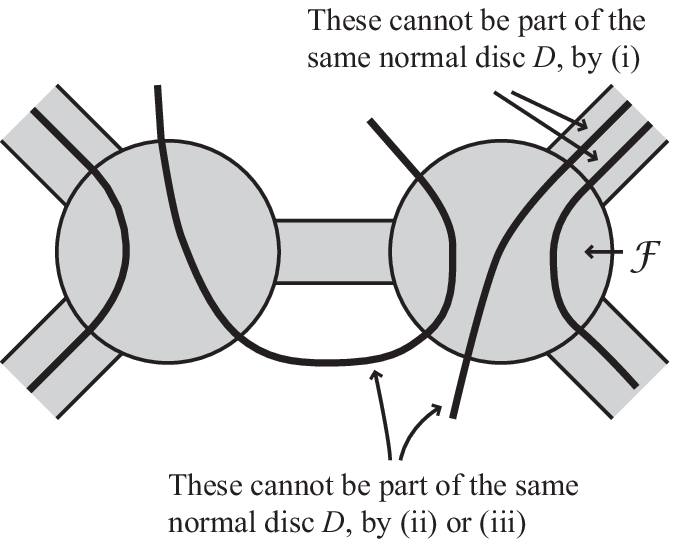}
}
\vskip 6pt
\centerline{Figure 7.}

When $T$ has no boundary and $H_0$ is a 0-handle such that $H_0 \cap {\cal F}$ is a thickening 
of the complete graph on 4 vertices, the above conditions imply that
$H_0 \cap T$ is a collection of {\sl triangles} and {\sl squares}, as shown in Figure 8.

\vskip 18pt
\centerline{
\epsfxsize=3in
\epsfbox{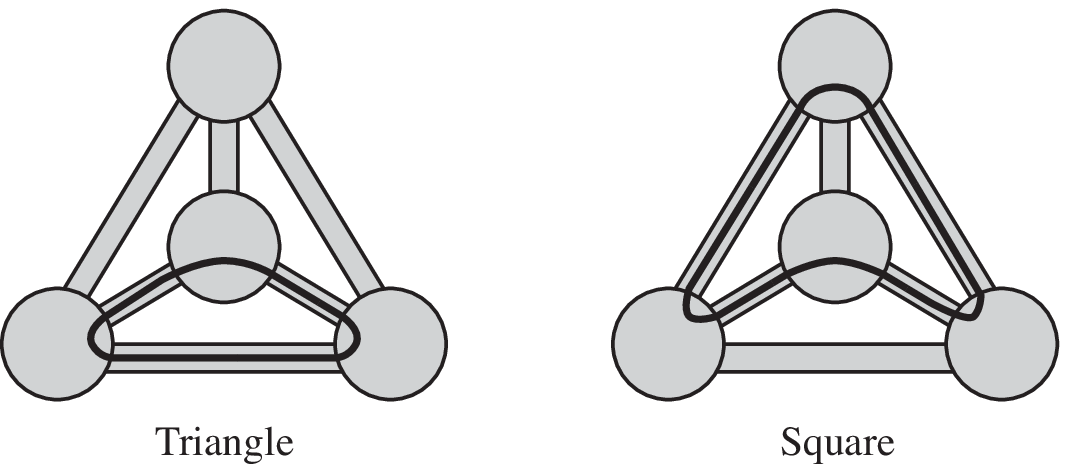}
}
\vskip 6pt
\centerline{Figure 8.}

We say that a simple closed curve properly embedded in $\partial X$ is {\sl standard} if
\item{(i)} it is disjoint from the 2-handles;
\item{(ii)} it intersects each 1-handle $D^1 \times D^2$ in $D^1 \times P$,
where $P$ is a finite set of points in $\partial D^2$;
\item{(iii)} it intersects ${\rm cl}(\partial {\cal H}^0_X - {\cal F})$ in
a collection of properly embedded arcs.

The following is Proposition 4.4 in [4]. It is a variant of a well-known result in normal
surface theory.

\noindent {\bf Proposition 3.4.} {\sl Let ${\cal H}_X$ be a handle structure on a compact
irreducible 3-manifold $X$. Let $T$ be a properly embedded, incompressible, boundary-incompressible
surface in $X$, with no 2-sphere components. Suppose that each component of $\partial T$ is
standard and intersects each component of $\partial X \cap {\cal H}_X^0$ and 
$\partial X \cap {\cal H}_X^1$ in at most one arc and no simple closed curves.
Then there is an ambient isotopy, supported in the interior of $X$, taking $T$ into
normal form.}

In the remainder of this section, $M$ will always be a compact orientable
3-manifold with a handle structure ${\cal H}_M$, and $S$ will be a compact subsurface of 
$\partial M$ such that $\partial S$ is standard. In this case, we say that ${\cal H}_M$ 
is a {\sl handle structure for the pair} $(M,S)$.

\noindent {\bf Definition 3.5.}
Let ${\cal H}_M$ be a handle structure for the pair $(M,S)$. 
A handle $H$ of ${\cal H}_M$ is a {\sl parallelity handle} if
it admits a product structure $D^2 \times I$ such that
\item{(i)} $D^2 \times \partial I = H \cap S$;
\item{(ii)} each component of ${\cal F}^0 \cap H$ and ${\cal F}^1 \cap H$
is $\beta \times I$, for a subset $\beta$ of $\partial D^2$.

We will typically view the product structure $D^2 \times I$
as an $I$-bundle over $D^2$.

The main example of a parallelity handle arises when $M$ is obtained by cutting a 3-manifold $X$
along a properly embedded, normal surface $T$, 
and where $S$ is the copies of $T$ in $M$. Then, if $T$ contains 
two normal discs in a handle that are normally 
parallel and adjacent, the space between them becomes a parallelity
handle in $(M,S)$.

The $I$-bundle structures on the parallelity handles can be chosen so that,
when two parallelity handles are incident, their $I$-bundle structures
coincide along their intersection. (This is established in the
proof of Lemma 5.3 in [4].) So, the union of the parallelity
handles forms an $I$-bundle over a surface $F$. This is termed the {\sl parallelity bundle} ${\cal B}$.
The $I$-bundle over $\partial F$ is termed
the {\sl vertical boundary} $\partial_v {\cal B}$ of ${\cal B}$,
and the $\partial I$-bundle over $F$ is called
the {\sl horizontal boundary} $\partial_h {\cal B}$.

As in [4], it will be technically convenient to consider enlargements of
such structures. These will still be an $I$-bundle over a surface
$F$, and near the $I$-bundle over $\partial F$, they will be a union
of parallelity handles, but elsewhere need not be. The precise definition
is as follows.

\noindent {\bf Definition 3.6.}
Let ${\cal H}_M$ be a handle structure for the pair $(M,S)$.
A {\sl generalised parallelity bundle} ${\cal B}$ is a 3-dimensional
submanifold of $M$ such that
\item{(i)} ${\cal B}$ is an $I$-bundle over 
a compact surface $F$;
\item{(ii)} the $\partial I$-bundle is ${\cal B} \cap S$;
\item{(iii)} ${\cal B}$ is a union of handles of ${\cal H}_M$;
\item{(iv)} any handle in ${\cal B}$ that intersects the $I$-bundle
over $\partial F$ is a parallelity handle, where the $I$-bundle
structure on the parallelity handle agrees with the $I$-bundle structure of
${\cal B}$;
\item{(v)} ${\rm cl}(M - {\cal B})$ inherits
a handle structure.

\noindent The $I$-bundle over $\partial F$ is termed
the {\sl vertical boundary} $\partial_v {\cal B}$ of ${\cal B}$,
and the $\partial I$-bundle over $F$ is called
the {\sl horizontal boundary} $\partial_h {\cal B}$.

We will also need to make some modifications to handle structures,
as follows. 

\noindent {\bf Definition 3.7.}
Let $G$ be an annulus properly embedded in $M$, with boundary in $S$.
Suppose that there is an annulus $G'$ in $\partial M$ such that $\partial G = \partial G'$.
Suppose also that $G \cup G'$ bounds a 3-manifold $P$ such that
\item{(i)} either $P$ is a parallelity region between $G$ and $G'$, or $P$ lies
in a 3-ball $B$ such that $B \cap \partial M$ is a disc;
\item{(ii)} $P$ is a non-empty union of handles;
\item{(iii)} ${\rm cl}(M - P)$ inherits
a handle structure from ${\cal H}_M$;
\item{(iv)} any parallelity handle of ${\cal H}_M$ that intersects
$P$ lies in $P$;
\item{(v)} $G$ is a vertical boundary component of a generalised
parallelity bundle lying in $P$.

\noindent Removing the interiors of $P$
and $G'$ from $M$ is called an {\sl annular
simplification}. Note that the resulting 3-manifold
$M'$ is homeomorphic to $M$, even though $P$ may
be homeomorphic to the exterior of a non-trivial knot
when it lies in a 3-ball. (See Figure 9.)
The boundary of $M'$ inherits a copy of $S$, which we denote
by $S'$, by setting $S' = (S \cap \partial M') \cup (\partial M' - \partial M)$.
Thus, $(M',S')$ is homeomorphic to $(M,S)$. Moreover, when $M$ is embedded 
within a bigger closed 3-manifold, then
$(M',S')$ is ambient isotopic to $(M,S)$.

\vskip 18pt
\centerline{
\epsfxsize=3in
\epsfbox{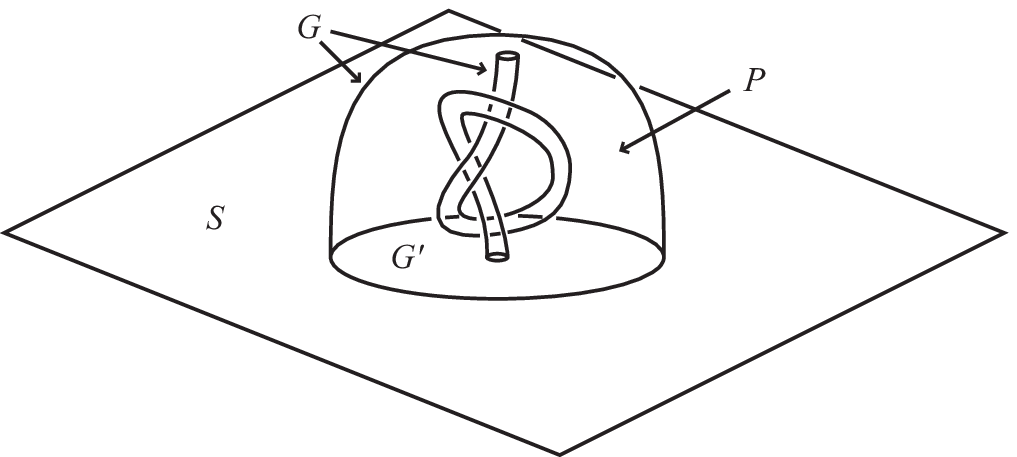}
}
\vskip 6pt
\centerline{Figure 9.}

This definition is similar to that given in [4]. However, there, $S$ was required to be
incompressible. We do not make that assumption here, but as a result,
condition (i) has been modified a little.

The following is Lemma 5.5 in [4].

\noindent {\bf Lemma 3.8.} {\sl Let ${\cal H}_M$ be a handle structure
for the pair $(M,S)$. Let ${\cal H}'_M$ be a handle structure obtained
from ${\cal H}_M$ by annular simplifications. Then any parallelity
handle for ${\cal H}_M$ that lies in ${\cal H}'_M$ is a parallelity
handle for ${\cal H}'_M$.}

The following is Corollary 5.7 in [4].

\noindent {\bf Theorem 3.9.} {\sl Let $M$ be a compact orientable irreducible
3-manifold with a handle structure ${\cal H}_M$. Let $S$ be an incompressible 
subsurface of $\partial M$ such that $\partial S$ is standard in $\partial M$. Suppose
that ${\cal H}_M$ admits no annular simplification. Then ${\cal H}_M$
contains a generalised parallelity bundle that contains every parallelity
handle and that has incompressible horizontal boundary.}

From this, we deduce the following.

\noindent {\bf Theorem 3.10.} {\sl Let $M$ be a compact connected orientable irreducible
3-manfiold, with boundary a collection of incompressible tori. Let $S$ be a
union of components of $\partial M$. Let ${\cal H}_M$ be a handle structure for $(M,S)$ that admits no
annular simplification. Suppose that there is no essential properly
embedded annulus in $M$, with boundary in $S$, and that lies entirely
in the parallelity handles. Suppose also that $M$ is not an $I$-bundle over a torus
or Klein bottle. Then $M$ contains a generalised parallelity
bundle ${\cal B}$ such that
\item{(i)} ${\cal B}$ contains every parallelity handle; and
\item{(ii)} ${\cal B}$ is a collection of $I$-bundles over discs.

}

This is very close to the statement of Proposition 5.8 in [4]. However,
there, $S$ was a collection of annuli, whereas here it is tori. In addition,
a slightly more precise hypothesis has been made about the annuli
properly embedded in $M$.

\noindent {\sl Proof.} By Theorem 3.9, there is a generalised parallelity
bundle ${\cal B}$ that contains every parallelity handle and 
with incompressible horizontal boundary. Let ${\cal B}'$
be the union of the components of ${\cal B}$ that are not
$I$-bundles over discs. Its horizontal boundary $\partial_h{\cal B}'$ is a
subsurface of $S$, and hence a collection of annuli and
tori. However, if there is any toral component of $\partial_h{\cal B}'$, then
this lies in a component of ${\cal B}'$ which is an
$I$-bundle over a Klein bottle or torus. This then
is all of $M$, which is contrary to hypothesis.
Hence, the horizontal boundary of ${\cal B}'$ just consists
of annuli. By hypothesis, each component of $\partial_v {\cal B}'$
is inessential. But $\partial_v {\cal B}'$ is incompressible,
and so each component of $\partial_v {\cal B}'$ is boundary parallel.
The proof now follows that of Proposition 5.8 in [4]
with almost no change. We note that when an annular
simplification is performed  at the end of the proof of
Proposition 5.8 in [4], an annular simplification can also be
done in our case. This completes the proof of the theorem.
$\square$

\vskip 18pt
\centerline{\caps 4. Proof of the main theorem}
\vskip 6pt

Let $K$ be a satellite knot with companion knot $L$. Let $X$ be the exterior of
$K$. 

\noindent {\sl Claim.} There is a companion knot $L'$ for $K$ with the following properties:
\item{(i)} $L'$ has $L$ as a connected summand (possibly $L' = L$);
\item{(ii)} $L'$ is not a non-trivial connected summand for any other companion for
$K$.

Suppose that $L$ fails condition (ii). In other words, suppose that
there is another companion $L_1$ for $K$, such that $L_1$ is the
connected sum of $L$ and some other non-trivial knot. Then $L_1$ is a
satellite of $L$, and so we have an inclusion $K \subset N(L_1) \subset N(L)$.
(See Figure 10.)
This knot $L_1$ may also fail condition (ii), but in this case we get
another knot $L_2$ such that $K \subset N(L_2) \subset N(L_1) \subset N(L)$,
and where $L_2$ has $L$ as a connected summand, and so on. Each torus $\partial N(L_i)$
is essential in $X$, and they are disjoint and non-parallel. By Kneser's
theorem, there is an upper bound on the number of such tori. Hence, eventually,
we obtain the required knot $L'$, as claimed.

\vskip 18pt
\centerline{
\epsfxsize=3.2in
\epsfbox{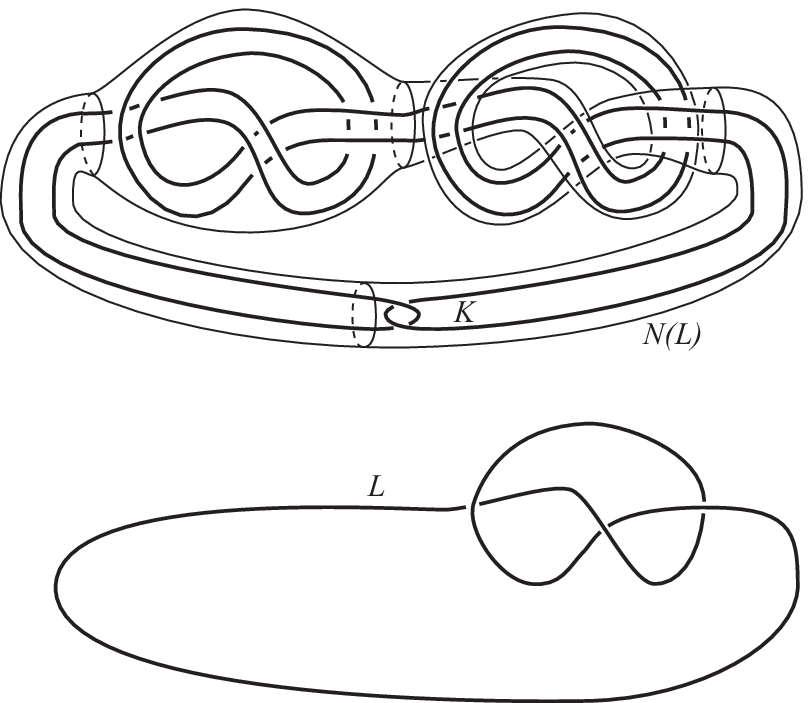}
}
\vskip 6pt
\centerline{Figure 10.}

Let $D$ be a diagram of $K$ with minimal crossing number. Our aim is to construct
a diagram $D'$ for $L'$ with
crossing number at most $3 \times 10^{10} \ c(D)$. Then, applying Theorem 1.2, we deduce that
$$c(L) \leq 152 \ c(L') \leq 152 \ c(D') \leq 152 \times 3 \times 10^{10} \ c(D) < 10^{13} \  c(K),$$
thereby proving Theorem 1.1. This will also prove Theorem 1.4, because if there
is no essential torus in $N(L) - {\rm int}(N(K))$, then we must have $L' = L$.

We may assume that $K$ does not have $L'$ as a connected summand,
because by Theorem 1.2, this would imply that $c(K) \geq c(L')/152$,
which is stronger than the required inequality.

Give $X$ the handle structure described in Section 2, which
we denote by ${\cal H}_X$. Let $V$ be the solid torus $N(L')$, and let
$T$ be the torus $\partial V$.
Since $T$ is incompressible in $X$, it may be placed in normal form with respect to ${\cal H}_X$.
Note that $T$ then inherits a handle structure where the $i$-handles
are $T \cap {\cal H}_X^i$.
Cutting $X$ along $T$ gives two 3-manifolds, one of which is a copy of the exterior
of $L'$, the other of which is $V$ with the interior of a small regular
neighbourhood of $K$ removed. Let $Y$ be this latter manifold. 
It inherits a handle structure ${\cal H}_Y$.

We wish to place an affine structure upon ${\cal H}_Y$. To do this, we start by 
straightening $T$ as much as possible, in the following way. In each 1-handle $D^1 \times D^2$
of ${\cal H}_X$, we realise $T \cap (D^1 \times D^2)$ as $D^1 \times \alpha$,
where $\alpha$ is a collection of straight arcs in the polygonal structure on $D^2$.
Then we make each arc of $T \cap {\cal H}^0_X \cap {\cal H}^2_X$ straight.
Thus the boundary of each normal disc of $T$ in each 0-handle is now a
concatenation of straight arcs. We then realise each normal disc of 
$T$ in ${\cal H}^0_X$ as a union of flat polygons in that 0-handle, as follows.

Observe first that although the polyhedron $P$ associated with each 0-handle
is not convex, it is star-shaped. In other words, there is a point $v$
in the interior of $P$, such that for each point on $\partial P$, 
the straight line joining it to $v$ lies within $P$. 
Moreover, the interior of this straight line lies within
the interior of $P$. We may therefore use a dilation about $v$,
with any positive scale factor less than $1$, to create a copy of $\partial P$
lying within the interior of $P$. Create a nested collection of such copies
of $\partial P$, the number of copies being equal to the number of
normal discs that we need to insert into $P$. We now create these normal
discs one at a time, starting with one that has boundary that is innermost
in $\partial P$. The boundary of this disc has already been specified
as some curve $C$ in $\partial P$. Use the star-shaped nature of $P$
to create an annulus interpolating between $C$ and a curve
on the outermost copy of $\partial P$. Since $C$ is a union of straight arcs,
this annulus is a union of flat quadrilaterals. The boundary
component of this annulus that lies on the dilated copy of $\partial P$
bounds a disc in this copy of $\partial P$ which is a union of 
flat polygons. The union of these with the annulus is the required
normal disc bounded by $C$. We now repeat this procedure for the
remaining curves of $T \cap \partial P$, starting with the next
innermost curve, and using the next copy of $\partial P$. In this way,
each normal disc of $T$ has been realised as a union of flat polygons.

We now want to bound the number of flat polygons that comprise each
normal disc. The boundary of each normal triangle is $6$ straight
arcs, and these give rise to $6$ flat polygons in the annular part
of the normal disc, together with $10$ flat polygons in the
copy of $\partial P$. This is $16$ in total. Similarly,
one can compute that each normal square is composed of $25$ flat
polygons.

Having realised $T$ in this way, we cut ${\cal H}_X$ along $T$,
and then ${\cal H}_Y$ inherits an affine handle structure. 
We view this as a handle structure for the pair $(Y,T)$.
Thus, whenever we consider a parallelity handle in ${\cal H}_Y$, its horizontal boundary
lies in the copy of $T$ in $Y$.

Each 0-handle $H_0$ of ${\cal H}_X$
gives rise to at most $6$ handles of ${\cal H}_Y$ that are not
parallelity handles. (This number $6$ arises when we cut the 0-handle
along $4$ normal triangles of distinct types and a normal square.)
Let $H_0'$ be the non-parallelity 0-handles of ${\cal H}^0_Y$ lying in $H_0$.
We wish to bound the number of faces in the polyhedral structure
of $H_0'$. Clearly, the maximal number of faces occurs when $H_0$
is an unexceptional 0-handle that is cut along all $4$ triangle types and a square type, and so it is
this configuration that we will examine.

The faces of $H_0$ come in the following types: the intersections with
${\cal H}^1_X$ (of which there are $4$), the intersections with ${\cal H}^2_X$
(of which there are $6$) and the intersections with $\partial X$ (of
which there are $4$ components, forming a total of $16$ faces).
Each face in $H_0 \cap {\cal H}^1_X$ gives rise to $5$ 
faces of $H_0' \cap {\cal H}^1_Y$. The faces of $H_0 \cap {\cal H}^2_X$
become $22$ faces of $H_0' \cap {\cal H}^2_Y$. The faces of $H_0 \cap \partial X$
stay as faces of $H_0'$. Finally, the normal triangles and squares of $T$ give rise
to $(16 \times 8) + (25 \times 2) = 178$ faces of $H_0'$. So, in total,
the number of faces of $H_0'$ is at most $(4 \times 5) + 22 + (4 \times 4) + 178 =
236$.

Apply as many annular simplifications to ${\cal H}_Y$ as possible,
creating a handle structure ${\cal H}_{Y'}$ for an isotopic
copy of $(Y,T)$, which we call $(Y',T')$. By Theorem 3.9, ${\cal H}_{Y'}$
contains a generalised parallelity bundle ${\cal B}$ that contains every parallelity
handle of ${\cal H}_{Y'}$ and that has incompressible horizontal boundary.
Note that ${\cal H}_{Y'}$ inherits an affine handle structure, since it is 
a union of handles of ${\cal H}_Y$.

Since $\partial_h {\cal B}$ is an incompressible subsurface of $T'$, it
is a collection of discs and annuli. Thus, ${\cal B}$ is a collection
of $I$-bundles over discs, annuli and M\"obius bands.

The proof now divides into two cases.

\noindent {\sl Case 1.} Some component of ${\cal B}$ is an $I$-bundle
over a M\"obius band.

This component of ${\cal B}$ is a solid torus, and its horizontal
boundary is an annulus that winds twice around the solid torus $Y' \cup N(K)$.
Each boundary component of this annulus is therefore a 2-cable of $L'$.
Recall that a knot $C$ is a {\sl 2-cable} of a knot $L'$ if $C$ is the boundary
of an embedded M\"obius band, the core of which is a copy of $L'$.
The linking number between $C$ and this core curve is the {\sl twisting number}
of the 2-cable. 

Let $C$ be one of these 2-cables.
We now wish to find an upper bound for the crossing number of $C$. We will do
this by finding a bound on the number of crossings in the diagram of $C$
that is obtained by projecting $C$ vertically onto the horizontal plane.

Now, $C$ has the structure of a cell complex, where each 1-cell is a component
of intersection between a handle of $\partial_h {\cal B}$ and a handle
of $T$ that does not lie in $\partial_h {\cal B}$. Each such 1-cell either lies in the
boundary of a 0-handle of ${\cal H}_X$ or it lies in a 1-handle of ${\cal H}_X$
and respects its product structure. When we project $C$ vertically, the
images of the latter 1-cells do not intersect each other or the image of any other
1-cell. Thus, the only way that crossings of the diagram of $C$ arise is
from the 1-cells of $C$ that lie in ${\cal H}^0_X$. Each such 1-cell lies in
the boundary of a normal disc of $T$ that does not lie in ${\cal B}$. There
are at most 10 such normal discs, consisting of at most 2 squares
and at most 8 triangles. Each square of $T$ has in its boundary 8
1-cells, each of which is straight in the polyhedral structure on the 0-handle 
of ${\cal H}_X$. Each triangle of $T$ has in its boundary 6 straight arcs.
Thus, in each 0-handle of ${\cal H}_X$, $C$ is composed of at most 
64 straight arcs. We may perturb these arcs a little so that, when projected
to the horizontal plane, the images of any two arcs intersect in at most 
one point. So, the resulting diagram of $C$ at most $64 \times 64 /2$ crossings
arising from each 0-handle of ${\cal H}_X$. There are $(4c(D) + 2)$ \break $0$-handles of $X$, which is at 
most $(14/3)c(D)$. So, the number of crossings of $C$ is
at most $9558 \ c(D)$.

Almost exactly the same argument gives that the modulus of the twisting of this 2-cable is
at most $9558 \  c(D)$. For we may consider a curve $C'$ in $T$
parallel to $C$. Then, the twisting of the cable is half the linking number
of $C$ and $C'$. The modulus of this linking number is at most the number of crossings between
$C$ and $C'$ when they are vertically projected. So, we get the same upper bound of
$9558 \  c(D)$ for the modulus of the twisting number.

The reason for bounding the crossing number of $C$ and the twisting of the
2-cable is that these are key quantities in the following theorem, which will be proved 
in Section 5.

\noindent {\bf Theorem 5.1.} {\sl Let $C$ be a knot that is a 2-cable of a knot $L'$ with twisting number $t$.
Then, $c(L') \leq 119024 \  (c(C) + |t|)$.}

Applying this result to our situation gives that $c(L') \leq 119024 \times 2 \times 9558 \  c(D) < 
3 \times 10^9 \ c(D)$,
which is better than the required bound. This proves Theorems 1.1 and 1.4 in this case.

\noindent {\sl Case 2.} No component of ${\cal B}$ is an $I$-bundle
over a M\"obius band.

Our aim is to construct a handle structure ${\cal H}_{V'}$ for an isotopic copy $V'$ of 
the solid torus $V$.
The proof now divides according to whether the conditions of Theorem 3.10
are met by ${\cal H}_{Y'}$, viewed as a handle structure for $(Y', T')$.

\noindent {\sl Case 2A.} There is no essential annulus properly embedded in $Y'$ 
with boundary in $T'$ and that lies entirely within the union of the
parallelity handles of ${\cal H}_{Y'}$.

By Theorem 3.10,  we may assume that ${\cal B}$ 
is a collection of $I$-bundles over discs. We replace each component
of ${\cal B}$ by a 2-handle, forming a handle structure ${\cal H}'_{Y'}$.
This inherits an affine structure from ${\cal H}_{Y'}$, since according
to the definition, we do not need to identify 2-handles with polyhedra.

The next step is to add a 2-handle as a thickened meridian disc
for $K$. This 2-handle is attached along the exceptional
0-handles and 1-handles. Finally, a 3-handle is added
in the remainder of $N(K)$. The result is a handle structure
${\cal H}_{V'}$ for $V' = Y' \cup N(K)$. We can give it an affine structure
as follows. We simply declare that the polyhedra associated with each
unexceptional $0$-handle and $1$-handle is the same as that in ${\cal H}'_{Y'}$.
The exceptional 0-handles and 1-handles must be modified slightly
to ensure that the new 2-handles intersect each of
these polyhedra in convex polygons. Thus, it is easy to arrange that the conditions in the
definition of an affine handle structure are satisfied by ${\cal H}_{V'}$.

\noindent {\sl Case 2B.} There is some properly embedded essential
annulus $A$ in $Y'$ with boundary in $T'$, and that lies in
the union of the parallelity handles of ${\cal H}_{Y'}$.

Then $A$ is properly embedded in the solid torus
$V' = Y' \cup N(K)$. Now, any properly embedded annulus in a solid torus $V'$ satisfies
at least one of the following:
\item{(i)} it is boundary parallel,
\item{(ii)} both its boundary curves are meridians and the annulus forms a knotted
tube joining these, or
\item{(iii)} at least one of its boundary curves bounds a disc in $\partial V'$. 

We claim that (ii) and (iii) do not arise. For suppose that (ii) holds.
Then $A$ cuts $V'$ into two pieces, one of which
is homeomorphic to the exterior of a non-trivial knot $K'$, and the other of which is a solid torus $W$.
Now, $K$ must lie in $W$ because the former piece
lies within a 3-ball in $V'$. Also, $K$ does not lie in a 3-ball in $W$,
since $K$ would then lie in a 3-ball in $V'$. In addition, $K$ is not
a core curve of $W$, because this would imply that $K$ has $L'$ as a connected
summand, and we are assuming that this does not occur. Thus, $K$ is a non-trivial satellite
of $L' \sharp K'$, which is contrary to our assumption about $L'$.
Now suppose that (iii) holds. Then the disc in $\partial V'$ that
is bounded by one of the components of $\partial A$ becomes a compression
disc for $A$ in $Y'$, which contradicts the assumption that $A$ 
is essential in $Y'$. This proves the claim.

Thus, $A$ is boundary parallel in $V'$. Since $A$ is essential in $Y'$,
$K$ must lie in the parallelity region $P$ between $A$ and a sub-annulus of $T'$.
We take $A$ to be innermost in $V'$, in the sense that any other essential
annulus properly embedded in $Y'$ with boundary in $T'$ and which lies within the union of
the parallelity handles must lie in $P$. Then $A$ is a vertical boundary component of ${\cal B}$.
Moreover, the component of ${\cal B}$ that is incident to $A$ lies in $P$.
This is because this component of ${\cal B}$ cannot be an $I$-bundle over
a M\"obius band, and so it is an $I$-bundle over an annulus.
So, if this component of ${\cal B}$ did not lie in $P$, this would
contradict the fact that we have taken $A$ to be innermost in $V$.

Remove $P - A$ from $V'$. The resulting 3-manifold
is an isotopic copy of $V'$, which we still call $V'$. It inherits an affine handle structure.
The intersection ${\cal B} \cap V'$ is a generalised parallelity bundle,
consisting of $I$-bundles over discs. 
Replace each of these $I$-bundles with a 2-handle, forming the
required affine handle structure ${\cal H}_{V'}$ for $V'$. Its only parallelity
handles are 2-handles.

Thus, in each of Cases 2A and 2B, we have created an isotopic copy $V'$
of $V$, and it has an affine handle structure ${\cal H}_{V'}$.

\vfill\eject
We claim that ${\cal H}_{V'}$ satisfies the conditions of Theorem 1.3,
apart from the fact that $X$ is not a solid torus.
Note that ${\cal H}_X$ satisfies these conditions. Then ${\cal H}_{V'}$ was
obtained from ${\cal H}_X$ by the following operations:
\item{(i)} cutting along a closed properly embedded normal surface;
\item{(ii)} removal of handles (but maintaining a handle structure);
\item{(iii)} replacing generalised parallelity bundles by 2-handles;
\item{(iv)} adding some handles (when filling in $N(K)$ in Case 2A).

It is clear that the conditions of Theorem 1.3 are preserved under (i), (ii) and (iii).
The addition of handles in (iv) also does not violate the conditions
of Theorem 1.3, because the 2-handles that are added are attached to the remnants of the 
exceptional 0-handles and 1-handles. The exceptional 1-handles of ${\cal H}_X$
each intersect intersect ${\cal H}^2_X$ in two components, and so the addition
of a further 2-handle does not violate the hypotheses of Theorem 1.3.
This proves the claim.

Each 0-handle and 1-handle
of ${\cal H}_{V'}$ is a handle of ${\cal H}_Y$. Thus, each 0-handle of
${\cal H}_{V'}$ lies in a 0-handle of ${\cal H}_X$. Moreover, each
0-handle of ${\cal H}_X$ contains at most six 0-handles of ${\cal H}_{V'}$.
This is because each 0-handle of ${\cal H}_X$ can support at most
5 types of triangles and squares of $T$ that are simultaneously disjoint.
These therefore divide the handle into at most six 0-handles that are not
parallelity handles.
The remaining 0-handles of ${\cal H}_Y$ are parallelity 0-handles,
and are therefore removed in the construction of ${\cal H}_{V'}$.

We now apply Theorem 1.3 to deduce that $V'$ has a core curve $C$ that
lies in the 0-handles and 1-handles of ${\cal H}_{V'}$, that
respects the product structure on the 1-handles, and that
intersects each 0-handle of ${\cal H}_{V'}$ in at most $48$
arcs. Moreover, these arcs are parallel to arcs $\alpha$ in the boundary
of the 0-handle, and each component of $\alpha$ intersects each face of
the 0-handle in at most $6$ straight arcs. We calculated above
that each 0-handle of ${\cal H}_X$ gave rise to at most six non-parallelity 0-handles
of ${\cal H}_Y$ and hence at most six 0-handles of ${\cal H}_{V'}$.
We also calculated that these $0$-handles of ${\cal H}_Y$ have at most $236$ faces.
Hence, the same is true of the 0-handles of ${\cal H}_{V'}$. (The exceptional
0-handles do not exceed this bound.) So, within each $0$-handle of ${\cal H}_X$, $C$ consists of at
most $6 \times 48 \times 236 = 67968$ straight arcs. 

We now project $C$ vertically to form a diagram $D'$ of $L'$. Now, $C$ lies
in the 0-handles and 1-handles of ${\cal H}_X$, and within the
1-handles, it respects their product structure. The interiors
of the 0-handles and the 1-handles have disjoint images under
the projection map. So, the only place that crossings of $D'$
can occur is within the images of the 0-handles of ${\cal H}_X$.
We have calculated that within each 0-handle, $C$ consists
of at most $67968$ straight arcs. These give rise to at most $(67968)^2$ crossings.
There are $(4c(D) + 2)$ \break $0$-handles of $X$, which is at most $(14/3)c(D)$.
So, the number of crossings of $D'$ is at most $(14/3)(67968)^2 c(D) < 3 \times 10^{10} c(D)$,
as required. 

Thus, in order to complete the proof of Theorems 1.1 and 1.4, all that remains
for us to do is prove Theorem 5.1. This we do in the next section.

\vskip 18pt
\centerline{\caps 5. The crossing number of 2-cables}
\vskip 6pt

In this section, we prove a version of the main theorem in a very special case,
where the satellite knot $K$ is a 2-cable. 

\noindent {\bf Theorem 5.1.} {\sl Let $K$ be a knot that is a 2-cable of a knot $L$
with twisting number $t$. Then, $c(L) \leq 119024 \  (c(K) + |t|)$.}

We may clearly assume that $L$ is a non-trivial knot, as otherwise, the statement
of the theorem is empty.

Let $X$ be the exterior of $K$.
Since $K$ is a 2-cable, it forms the boundary of an embedded M\"obius band, the core of which
is copy of $L$. Let $F_0$ be the restriction of this M\"obius band to $X$.

Let $D$ be a diagram for $K$ with minimal crossing number.
We now modify $D$ by performing Type I Reidemeister moves
which introduce kinks. We perform enough of these moves so
that the writhe of the new diagram $D'$ is equal to twice the twisting of the
2-cable. Thus, $D'$ has crossing number at most $2c(K) + 2 |t|$. Using this diagram,
give $X$ the affine handle structure described in Section 2, but without the
introduction of the exceptional handles. We denote this by ${\cal H}_X$.
Note that the number of 0-handles of ${\cal H}_X$ is $4c(D')$, which is
at most $8 (c(K) + |t|) $.

We may pick a simple closed curve on $\partial N(K)$ that winds once along $N(K)$
and that has blackboard framing with respect to $D'$. We may arrange that it is
standard in the handle structure and intersects each 0-handle of ${\cal H}_X$ in at most two arcs.
Moreover, it intersects each component of ${\cal H}^0_X \cap \partial X$
and ${\cal H}^1_X \cap \partial X$ in at most one arc.
Since the writhe of $D'$ is equal to twice the twisting of the 2-cable, this simple
closed curve is ambient isotopic to $\partial F_0$. Thus, we may arrange $F_0$ so that its
boundary is equal to this curve.

Because $L$ is non-trivial, $F_0$ is boundary-incompressible, and it is
also incompressible. So, by Proposition 3.4, there is an ambient isotopy,
supported in the interior of $X$, taking $F_0$ to a normal surface.
This does not move the boundary of $F_0$. 

Now let $F$ be a M\"obius band properly embedded in $X$, such that
\item{(i)} $F$ is normal;
\item{(ii)} $\partial F = \partial F_0$;
\item{(iii)} $F$ is ambient isotopic to $F_0$.

Choose $F$ so that the pair $(|F \cap {\cal H}_X^2|, |F \cap {\cal H}^0_X|)$ is minimal
among all properly embedded M\"obius bands satisfying the above three conditions. 
Here, we are placing lexicographical ordering on such pairs.
Thus, we first minimise $|F \cap {\cal H}_X^2|$, and then once this
has smallest possible value, we minimise $|F \cap {\cal H}^0_X|$.

Let $N(F)$ be a thin regular neighbourhood of $F$ in $X$, which is an $I$-bundle over $F$,
in which $F$ lies as a zero section. Let $\tilde F$ be the associated $\partial I$-bundle over
$F$, which is the annulus ${\rm cl}(\partial N(F) - \partial X)$.

Since $F$ is normal, it inherits a handle structure, where the $i$-handles of $F$ are the
components of ${\cal H}_X^i \cap F$. Similarly, $\tilde F$ inherits
a handle structure.

Let $M$ be the result of cutting $X$ along $F$. Then $M$ inherits a handle structure
${\cal H}_M$.  Note that $\tilde F$ is a subsurface of $\partial M$ and its boundary is
standard in ${\cal H}_M$. Thus, ${\cal H}_M$ is a handle structure for the
pair $(M, \tilde F)$. Let ${\cal B}$ be its parallelity bundle.

Note that because $\partial F$ runs over each component of ${\cal H}_X^0 \cap \partial X$ 
and ${\cal H}^1_X \cap \partial X$ at most once,
no parallelity handle of ${\cal H}_{M}$ intersects $\partial \tilde F$.
Hence, $\partial_h {\cal B}$ lies in the interior of $\tilde F$.
Let $\tilde \Gamma$ be the boundary of $\partial_h {\cal B}$, which is therefore a
collection of simple closed curves in $\tilde F$. 
We give $\tilde \Gamma$ a cell structure, where each 1-cell is a component of intersection between
adjacent handles of $\tilde F$. Let $\Gamma$ be the
image of $\tilde \Gamma$ in $F$ under the bundle map $\tilde F \rightarrow F$. 
Then $\Gamma$ is also a 1-complex. 

\noindent {\sl Case 1.} $\Gamma$ contains a core curve $C$ of $F$ as a subcomplex.

Note that $\Gamma$ has controlled intersection with each handle of ${\cal H}_X$,
in the following sense.  Each 1-cell of $\Gamma$ lies in the image of a 1-cell of $\tilde \Gamma$.
This is a component of intersection between two handles of $\tilde F$. One, $H_1$ say, lies in $\partial_h {\cal B}$,
and the other, $H_2$, does not lie in $\partial_h {\cal B}$. We will now
control the possibilities for $H_2$ within each 0-handle of ${\cal H}_X$.
Now in each 0-handle of ${\cal H}_X$, at most two normal discs of $F$ intersect
$\partial X$. Suppose that there is such a disc $E$. Then each 1-cell of
$E \cap \Gamma$ arises as a component of $E \cap {\cal F}^0$ or $E \cap {\cal F}^1$.
We may arrange that each such 1-cell is straight in the affine structure on the 0-handle.
Now, $E$ runs over each component of ${\cal F}^1$ at most once, and so
this gives rise to at most 6 1-cells of $E \cap \Gamma$. Between these
we have at most $5$ components of $E \cap {\cal F}^0$ which miss $\partial X$.
So, this gives at most $11$ 1-cells of $\Gamma$ arising from $E$.
The remaining normal discs of $F$ are squares and triangles, and so 
at most $10$ of these can give a handle $H_2$ of $\tilde F$ not in $\partial_h {\cal B}$.
At most $2$ of these are squares, and at most $8$ are triangles.
So, this gives at most $(2 \times 8) + (8 \times 6) = 64$ 1-cells of
$\Gamma$. Again, we may arrange that each of these is straight.
So, in each 0-handle of ${\cal H}_X$, we have at most $86$ straight 1-cells of
$\Gamma$. Projecting the core curve $C$, we obtain a 
diagram for $L$. The crossings of this diagram occur only in the
projections of the 0-handles of ${\cal H}_X$. There are at most $8(c(K) + |t|)$
such 0-handles. So, we get a diagram for $L$ with at most ${1 \over 2} \times (86 \times 85) \times 8 (c(K) + |t|)
= 29240 (c(K) + |t|)$ crossings, which is better
than the required bound.

\noindent {\sl Case 2.} $\Gamma$ does not contain a core curve of $F$ as a subcomplex.

We claim that we may find a core curve $C$ of $F$ which avoids $\Gamma$. To see this,
pick an ordering on the 1-cells of $\Gamma$, and remove an open regular neighbourhood
of these 1-cells from $F$ one at a
time. At each stage, we examine the complementary regions. Initially, this is just
a M\"obius band. We will show that, at each stage, one complementary region is
a M\"obius band with a (possibly empty) collection of open discs removed. 
We call this a punctured M\"obius band. Moreover, a core curve of this punctured
M\"obius band is also a core curve of $F$. Note this complementary region
either contains all of $\partial F$ or is disjoint from $\partial F$ because
$\Gamma$ is disjoint from $\partial F$.
As each new 1-cell of $\Gamma$ is removed from $F$, there are three options. It may
be completely disjoint from the previous cells, in which case this just punctures
one of the complementary regions. It may have just one endpoint incident to
the previous cells, in which case the complementary regions are unchanged up
to ambient isotopy. The main case is where both endpoints of the 1-cell
are incident to earlier cells. In this case, we cut a complementary region
along a properly embedded arc $\alpha$. We are only interested in the case where
this region is the punctured M\"obius band. If this joins different
punctures, the result is still a punctured M\"obius band. 
If $\alpha$ joins a puncture to the boundary of the punctured
M\"obius band, then cutting along $\alpha$ still results in
a punctured M\"obius band. If $\alpha$ joins the boundary to itself,
the result is an arc properly embedded in the M\"obius band. If the
arc is essential, then we have a core curve of $F$ as a subcomplex of $\Gamma$, which
is contrary to hypothesis. So, the arc is inessential, and there remains a complementary
region that is a punctured M\"obius band. If the arc joins a puncture to itself,
the result is a simple closed curve embedded within the M\"obius band.
This is either inessential, in which case there remains a complementary region that
is a punctured M\"obius band, or it is essential, in which case this leads to a core curve
of $F$ as a subcomplex of $\Gamma$,
and again this is contrary to hypothesis. This proves the claim.

Let $C$ be this core curve. After an isotopy in the complement of $\Gamma$, we may assume that it intersects
each handle of $F$ in at most one properly embedded arc.
Let $M_1$ be the union of the $I$-fibres in $N(F)$ over $C$. This is a M\"obius band.
Let $\tilde C$ be the boundary of $M_1$. Then $\tilde C$ misses $\tilde \Gamma$.
It therefore misses $\partial_h {\cal B}$ or lies in the interior of $\partial_h {\cal B}$.
The argument divides into these two cases.

\noindent {\sl Case 2A.} $\tilde C$ misses $\partial_h {\cal B}$.

Then $C$ lies in the surface $F_-$ that is obtained from $F$ by removing
the interior of the image of $\partial_h {\cal B}$ under the bundle map $\tilde F \rightarrow F$.
Note that whenever an $i$-handle of $\tilde F$ lies in $\partial_h {\cal B}$, 
so do all the $j$-handles with $j > i$ that are incident to it.
Hence, $F_-$ inherits a handle structure. We may therefore ensure that
$C$ misses the 2-handles of $F_-$ and respects the product structure on
the 1-handles. We may also isotope $C$ in $F_-$
so that it intersects each handle of $F_-$ in at most one properly embedded arc. 
Because $\tilde C$ misses $\partial_h {\cal B}$, 
each handle of $F$ that intersects $C$ is disjoint from ${\cal B}$ on
both sides. There can be at most 7 such normal discs of $F$ in each 0-handle of $X$:
at most 4 triangles, at most one square and at most two further normal discs that intersect
$\partial X$. As in Section 4, we may arrange that each triangle is made up 16 flat
polygons and that the square is made up of 25 flat polygons. Similarly, the normal
discs that intersect $\partial X$ consist of at most 84 polygons in total. We may ensure that
$C$ intersects each of the polygons in at most one straight arc. So, in each 0-handle of ${\cal H}_X$,
$C$ is at most $(4 \times 16) + 25 + 84 = 173$ straight arcs.
Therefore, the projection of $C$ has at most ${1 \over 2} (173 \times 172) \times 4 \ c(D') \leq
119024 \ (c(K)+|t|)$ crossings, as required.

\noindent {\sl Case 2B.} $\tilde C$ lies in $\partial_h {\cal B}$.

Our aim here is to reach a contradiction. Let ${\cal B}'$ be the component of
${\cal B}$ that contains $\tilde C$. Then $\partial_h {\cal B}'$ is
either connected or disconnected, depending of whether the base surface of
the $I$-bundle ${\cal B}'$ is non-orientable or orientable.
We consider these two cases separately.

\noindent {\sl Case 2B(i).} $\partial_h{\cal B}'$ is connected.

Then the base surface of ${\cal B}'$ is non-orientable, and therefore
contains a properly embedded, orientation-reversing simple closed curve. The union
of the fibres in ${\cal B}'$ over this curve is a M\"obius band $M_2$.
Its boundary is a simple closed curve in the annulus $\tilde F$. 
This cannot bound a disc in $\tilde F$, for then the union of
$M_2$ with this disc would be an embedded projective plane
in $S^3$, which is well known to be impossible. Thus,
$\partial M_2$ is a core curve of $\tilde F$ and therefore
separates $\tilde F$ into two annuli. Attach one of these
annuli to $M_2$ to form a properly embedded M\"obius band $M_3$
in $X$. The boundary curves of $M_3$ and $F$ are disjoint and
therefore cobound an annulus in $\partial X$. Attach this annulus
to $M_3 \cup F$ to form an embedded Klein bottle in $S^3$. This again
is a contradiction.

\noindent {\sl Case 2B(ii).} $\partial_h {\cal B}'$ is disconnected.

Then, ${\cal B}'$ is a product $I$-bundle. Let $A_1$ be the union of the $I$-fibres in ${\cal B}'$ that are 
incident to $\tilde C$. This is an annulus. The boundary curve
$\partial A_1 - \tilde C$ is a simple closed curve in $\tilde F$. It cannot bound
a disc in $\tilde F$, for then the union of this disc with $A_1 \cup M_1$
would be an embedded projective plane in $S^3$. So, $\partial A_1 - \tilde C$ divides $\tilde F$
into two annuli. Let $A_2$ be the annulus that does not contain $\tilde C$.
Let $A_3$ be the sub-annulus of $\tilde F$ lying between the
two components of $\partial A_1$.

Note that $\tilde F - \tilde C$ consists of two annuli. The restriction of the
bundle map $\tilde F \rightarrow F$ to each of these annuli is an injection.
In particular, the restriction of $\tilde F \rightarrow F$ to $A_2 \cup A_3$ is an injection.

We claim that $\partial A_1 - \tilde C$ intersects each handle of $\tilde F$
in at most one properly embedded arc. Now, when a handle of $\tilde F$ intersects $\tilde C$, 
it does so in one properly embedded arc. Since $\tilde C$ lies in ${\cal B}$, this handle
of $\tilde F$ runs parallel to another handle of $\tilde F$, with the region
between these two handles lying in ${\cal B}$. So, $\partial A_1 - \tilde C$
intersects this other handle of $\tilde F$ in just one arc. 

Thus, each handle of $\tilde F$ is divided into at most two components by
$\partial A_1 - \tilde C$ and at most one of these lies in $A_2$. So,
each $i$-handle of $F$ gives rise to at most one component of
$A_2 \cap {\cal H}^i_X$. Hence, for $i=0,1$ and $2$,
$|A_2 \cap {\cal H}^i_X| \leq |F \cap {\cal H}^i_X|$.

We claim that this inequality is strict for $i=0$.
Note that whenever an $i$-handle of $\tilde F$ lies in ${\cal B}$, so does
every $j$-handle adjacent to it, provided $j > i$. So, ${\rm cl}(\tilde F - \partial_h {\cal B})$
has a handle structure. In particular, each component of ${\rm cl}(\tilde F - \partial_h {\cal B})$
contains a 0-handle of $\tilde F$. So, $A_3$ contains at least
one 0-handle of $\tilde F$. Hence, $|A_2 \cap {\cal H}^0_X| < |F \cap {\cal H}^0_X|$,
as claimed.
 
Let $F'_0$ be the M\"obius band $M_1 \cup A_1 \cup A_2$.
This is properly embedded in $X$, and after a small isotopy,
we may arrange that it has the same boundary as $F$.
It has strictly fewer components of intersection
with the 0-handles of $X$ than $F$, and at most as many components
of intersection with ${\cal H}^2_X$. We have verified
this for $A_2$. But $M_1$ is composed of $I$-fibres in $N(F)$,
and $A_1$ is composed of $I$-fibres in ${\cal B}$. Hence,
adding these to $A_2$ does not increase the number of
components of intersection with ${\cal H}_X^0$ or ${\cal H}_X^2$. Now, $F_0'$
might not be a normal surface, but when we apply the usual normalisation procedure
(leaving the boundary fixed), the complexity of the pair
$(|F_0' \cap {\cal H}_X^2|, |F_0' \cap {\cal H}_X^0|)$ does not increase.
Thus, we end with a normal M\"obius band properly embedded in $X$, 
with the same boundary as $F$, but with 
smaller complexity. This is ambient
isotopic to $F$ by the following result. Hence this gives a contradiction.

\noindent {\bf Lemma 5.2.} {\sl Let $K$ be a knot in the 3-sphere. Let $F$
and $F_0$ be essential M\"obius bands properly embedded in the exterior
of $K$ with equal boundaries. Then $F$ and $F_0$ are
ambient isotopic in the exterior of $K$.}

\noindent {\sl Proof.} We make use of the JSJ decomposition of the exterior $X$ of $K$. The JSJ tori
are a collection of disjoint properly embedded essential tori, with the property that
any other properly embedded essential torus can be ambient isotoped off them.
The JSJ tori divide $X$ into pieces that are Seifert fibred or atoroidal. Since $X$ is a knot
exterior, the possibilities for the Seifert fibred pieces are very limited (see for example [1]).
Each has base surface that is a planar surface. Moreover, if the Seifert fibered
piece has any singular fibres, the base orbifold is either a disc with two singularities
of coprime order, or an annulus with one singularity.

Let $W$ be a regular neighbourhood of $F \cup \partial N(K)$ in $X$. This is
a Seifert fibre space. The boundary component $\partial W - \partial N(K)$ 
either bounds a solid torus with interior disjoint from $W$ or is essential in $X$.
In both cases, it can be
ambient isotoped off the JSJ tori, and then $W$ lies in a Seifert fibred piece $R$ of the
JSJ decomposition.
We can then arrange that the $W$ is a union of fibres in $R$, and that
the M\"obius band $F$ is also a union of fibres. Thus, $F$ contains a singular
fibre of order 2, and it projects to an embedded arc in the base orbifold of $R$.
This arc runs from $\partial R$ to the order 2 singularity. Note that there is just one
order 2 singularity in this orbifold, because if there is another, it has coprime order.
Now, removing the singularities from the orbifold gives a pair of pants $P$, and
there is, up to ambient isotopy, a unique properly embedded arc
in $P$ joining any two boundary components.

Now let $W'$ be a regular neighbourhood of $F_0 \cup \partial N(K)$ in $X$.
Then, $W'$ also lies in $R$ after an ambient isotopy, and it too is a union of
fibres in $R$, as is $F_0$. Thus, after removing the singular fibres,
$F_0$ projects to an arc in $P$ that is isotopic to the previous one.
Hence, $F_0$ is ambient isotopic to $F$. $\square$

This completes the proof of Theorem 5.1, and hence Theorems 1.1 and 1.4. $\square$

\vskip 18pt
\centerline{\caps 6. A related result}
\vskip 6pt

Recall that a knot $\tilde L$ is a {\sl cable} of a knot $L$ if $\tilde L$ is a
simple closed curve on $\partial N(L)$ that does not bound a disc
in $N(L)$.

\noindent {\bf Theorem 1.5.} {\sl Let $K$ be a satellite knot with companion knot $L$.
Suppose that $L$ is prime. Then, for some knot $\tilde L$, which is either $L$ or a cable of $L$,
$$c(K) \geq {c(\tilde L) \over 152}.$$}

\noindent {\sl Proof.} Let $D$ be a diagram for $K$ with minimal crossing number.
Let $X$ be the exterior of $K$.
Give $X$ the handle structure described in Section 2, but without
the introduction of the exceptional handles. Call this handle structure
${\cal H}_X'$. We do not give it an affine structure, but instead, realise
the handles as shown in Figures 2 and 3 and the left of Figure 5.
Let $T$ be the torus $\partial N(L)$. This may be placed in normal form
with respect to ${\cal H}'_X$. Cutting $X$ along $T$ gives two 3-manifolds,
one of which is a copy of the exterior of $L$. Call this latter manifold $Z$,
and let ${\cal H}_Z$ be the handle structure that it inherits.
View this as a handle structure for the pair $(Z,T)$.

Apply as many annular simplifications to ${\cal H}_Z$ as possible,
creating a handle structure ${\cal H}_Z'$. According to Theorem 3.9,
${\cal H}_Z'$ has a generalised parallelity bundle ${\cal B}$ that
contains every parallelity handle of ${\cal H}_Z'$ and which has
incompressible horizontal boundary $\partial_h {\cal B}$.

Now $\partial_h {\cal B}$ cannot be all of $\partial Z$. Because $Z$
would then be an $I$-bundle over a torus or Klein bottle. Hence,
$\partial_h {\cal B}$ is a (possibly empty) collection of discs
and a (possibly empty) collection of annuli.

\noindent {\sl Claim.} If $\partial_h {\cal B}$ contains any
annuli, then these do not have meridional slope.

Suppose that, on the contrary, these annuli have meridional
slope. The vertical boundary of the corresponding components of
${\cal B}$ are then incompressible annuli properly embedded in $Z$,
with meridional boundary. They cannot be essential, because
$L$ would then be a composite knot. Hence, they are boundary parallel.
Note that, in this case, no component of ${\cal B}$ can be an
$I$-bundle over  a M\"obius band. For then ${\cal B}$ would contain
a properly embedded M\"obius band that is a union of $I$-fibres,
and this could be capped off with discs to form an embedded projective
plane in ${\Bbb R}^3$, which is known not to exist. Thus,
the vertical boundary components of ${\cal B}$ with meridional boundary
come in pairs which lie in the same component of ${\cal B}$.
We may therefore pick one such annulus $A$, with the property that
the component of ${\cal B}$ that it lies in is part of the parallelity
region between $A$ and $\partial Z$. Removing this parallelity region
(apart from $A$ itself) is therefore an annular simplification that can be made to ${\cal H}_Z'$,
which is contrary to hypothesis. This proves the claim.

Let $F$ be the surface $\partial Z - {\rm int}(\partial_h {\cal B})$.
This is either a collection of punctured annuli or a punctured torus,
where the punctures arise from components of ${\cal B}$ that are
$I$-bundles over discs. It inherits a handle structure. In the case where 
it is a collection of punctured annuli,
we may pick a core curve $\alpha$ of one of these annuli, which misses the
punctures, which runs only
over the 0-handles and 1-handles, which respects the product structure
on the 1-handles, and which intersects each handle in at most one arc.
For we may start with a core curve $\alpha$ of one of the annuli, slide it
off the punctures and the 2-handles, and then straighten it in the 1-handles. If $\alpha$
intersects some handle in more than one arc, then we may find an
arc $\beta$ in the handle, joining distinct arcs of $\alpha$ in that
handle, and with interior disjoint from $\alpha$. One can then cut
$\alpha$ at the two points of $\partial \beta$, remove one of the 
resulting arcs, and replace it by $\beta$. The result is still
a core curve of the annulus, but which intersects fewer handles.

The other case is when $F$ is a punctured torus. This time we pick a curve $\alpha$ which avoids 
the punctures and the 2-handles,
and which respects the product structure on the 1-handles,
and which has longitudinal slope, say, on $\partial Z$. This time,
it may not be possible to isotope $\alpha$ so that it runs
over each handle in at most one arc. The modifications described
above will not necessarily keep $\alpha$ as a longitude.
However, they will not change its class in $H_1(\partial Z; {\Bbb Z}/2)$.
So, we may still find a simple closed curve $\alpha$ on $F$ which
runs over each handle in at most one arc, which misses the 2-handles and which
respects the product structure on the 1-handles, which is essential in $\partial Z$
and which is not a meridian.

Let $\tilde L$ be this curve $\alpha$. It is either isotopic to $L$
or is a cable of $L$. 

We now have to be a little more precise about
the position of $\tilde L$. We first arrange that each normal disc of
$T$ sits within the handle of ${\cal H}'_X$ that contains it as described in
Section 6.5 of [4]. Then we arrange that, whenever $\tilde L$
runs over one of these normal discs, then it does so as described
in Section 6.5 of [4]. We then project $\tilde L$ vertically,
forming a diagram $\tilde D$, and we have to bound its crossing number.
The details of the argument are identical to those in Section 6.6 of [4].
In particular, the argument there gives that $c(\tilde D) \leq 152 \ c(D)$.
This proves the theorem, because
$$c(\tilde L) \leq c(\tilde D) \leq 152 \ c(D) = 152 \ c(K). \eqno{\square} $$

We close with some final remarks about the nature of this proof.
The arguments behind Theorems 1.1 and 1.5 are similar, but in the
former case, we focused on $Y = N(L) - {\rm int}(N(K))$,
whereas in the latter case, we used $Z = S^3 - {\rm int}(N(L))$.
Although the use of $Z$ leads to better a better constant,
one loses track of where the meridian of $Z$ lies.
It therefore seems very hard to avoid
the possibility that the knot $\tilde L$ we are considering
might be a cable of $L$. Only by considering $Y$ and using
Theorem 1.3 does it seem feasible to bypass this issue of cabling.

Nonetheless, Theorem 1.5 reduces the general problem of finding a lower bound on 
the crossing number of a satellite knot to the same problem for cables.
It is conceivable that it can be used as part of an alternative proof of
Theorem 1.1.

\vskip 18pt
\centerline{\caps References}
\vskip 6pt

\item{1.} {\caps R. Budney}, {\sl JSJ-decompositions of knot and link complements in $S^3$.}
Enseign. Math. (2) 52 (2006), no. 3-4, 319--359.

\item{2.} {\caps M. Freedman, Z-X. He,} {\sl Divergence-free fields: energy and asymptotic crossing number.}
Ann. of Math. (2) 134 (1991), no. 1, 189--229.

\item{3.} {\caps R. Kirby,} {\sl Problems in low-dimensional topology},
 AMS/IP Stud. Adv. Math., 2.2,  Geometric topology (Athens, GA, 1993),  35--473, Amer. Math. Soc., Providence, RI, 1997.

\item{4.} {\caps M. Lackenby}, {\sl The crossing number of composite knots}, J. Topology 2 (2009) 747--768.

\item{5.} {\caps M. Lackenby}, {\sl Core curves of triangulated solid tori}, Preprint.

\item{6.} {\caps S. Matveev,} {\sl Algorithmic topology and classification of $3$-manifolds},
Algorithms and Computation in Mathematics, Volume 9, Springer (2003).

\vskip 12pt
\+ Mathematical Institute, University of Oxford, \cr
\+ 24-29 St Giles', Oxford OX1 3LB, United Kingdom. \cr

\end